\documentclass[11pt]{tran-l}
\usepackage[ansinew]{inputenc}
\usepackage[centertags]{amsmath}
\usepackage{exscale}
\usepackage{relsize}
\usepackage{amsfonts}
\usepackage{amssymb}
\usepackage{amsthm}
\usepackage{newlfont}

\mathchardef\ordinarycolon\mathcode`\:
  \mathcode`\:=\string"8000
  \begingroup \catcode`\:=\active
    \gdef:{\mathrel{\mathop\ordinarycolon}}
  \endgroup
\newtheorem*{mthm*}{Main Theorem}

\newtheorem{thm}{Theorem}

\newtheorem{prop}[thm]{Proposition}

\newtheorem{Cor}[thm]{Corollary}
\newtheorem*{problem*}{Problem}
\newtheorem{conj}{Conjecture}

\begin{document}
\title{Steiner $t$-designs for large $t$}

\author{Michael Huber}

\address{Institut f\"{u}r Mathematik, MA6-2, Technische Universit\"{a}t Berlin, Stra{\ss}e des 17. Juni~136,
D-10623 Berlin, Germany}

\email{mhuber@math.tu-berlin.de}




\date{July 11, 2008}

\commby{}


\begin{abstract}
One of the most central and long-standing open questions in
combinatorial design theory concerns the existence of Steiner
\mbox{$t$-designs} for large values of $t$. Although in his
classical 1987 paper, L.~Teirlinck has shown that non-trivial
\mbox{$t$-designs} exist for all values of $t$, no non-trivial
Steiner \mbox{$t$-design} with $t>5$ has been constructed until now.
Understandingly, the case $t=6$ has received considerable attention.
There has been recent progress concerning the existence of highly
symmetric Steiner \mbox{$6$-designs}: It is shown in [M. Huber,
\emph{J. Algebr. Comb.} \textbf{26} (2007), pp. 453--476] that no
non-trivial flag-transitive Steiner \mbox{$6$-design} can exist. In
this paper, we announce that essentially also no block-transitive
Steiner \mbox{$6$-design} can exist.
\end{abstract}

\maketitle

\section{Introduction}\label{intro}

One of the most central and long-standing open questions in
combinatorial design theory concerns the existence of Steiner
\mbox{$t$-designs} for large values of $t$. Although in his
classical 1987 paper, L.~Teirlinck~\cite{Teir1987} has shown that
non-trivial \mbox{$t$-designs} exist for all values of $t$, no
non-trivial Steiner \mbox{$t$-design} with $t>5$ has been
constructed until now. Understandingly, the case $t=6$ has received
considerable attention. There has been recent progress concerning
the existence of highly symmetric Steiner \mbox{$6$-designs}: The
author~\cite{hu07} showed that no non-trivial flag-transitive
Steiner \mbox{$6$-design} can exist. Moreover, he classified all
flag-transitive Steiner \mbox{$t$-designs} with $t>2$
(see~\cite{Hu2001,Hu2005,Hu_Habil2005,hu07,hu07a} and~\cite{Hu2008}
for a monograph). These results answer a series of 40-year-old
problems and generalize theorems of J.~Tits~\cite{Tits1964} and
H.~Lüneburg~\cite{Luene1965}. Earlier, F.~Buekenhout,
A.~Delandtsheer, J.~Doyen, P.~Kleidman, M.~Liebeck, and
J.~Saxl~\cite{Buek1990,Del2001,Kleid1990,Lieb1998,Saxl2002} had
essentially characterized all flag-transitive Steiner
\mbox{2-designs}. All these classification results rely on the
classification of the finite simple groups.

In this paper, we announce that essentially no block-transitive
Steiner \mbox{$6$-design} can exist. This confirms a far-reaching
conjecture of P.~Cameron and C.~Praeger~\cite{CamPrae1993}, stating
that there are no non-trivial block-transitive \mbox{$6$-designs},
for the important case of Steiner designs. Consequently, a further
significant step towards an answer to the fundamental open question
\emph{``Does there exist any non-trivial Steiner
\mbox{$6$-design}?''} is provided -- at least in the case of highly
symmetric designs, quoting arguably Gian-Carlo Rota~\cite{Rota1964}:
\begin{quote}
``A combinatorial object without symmetries doesn't exist - by
definition.''
\end{quote}


\section{Combinatorial Designs}\label{designs}

The study of combinatorial designs deals with a crucial problem of
combinatorial theory, that of arranging objects into patterns
according to specified rules. This is a subject of considerable
interest in discrete mathematics and computer science, amongst
others. In particular, there are close connections of design theory
with graph theory~\cite{cali91,ton88}, finite and incidence
geometry~\cite{buek_handb95,demb68}, group
theory~\cite{cam99a,carm37,dimo96,wiel64}, coding
theory~\cite{cali91,Hu_cod2008,hupl_handb98,hupl03},
cryptography~\cite{pei06,stin05}, as well as classification
algorithms~\cite{kaos06}.

Combinatorial designs may be regarded as generalizations of finite
projective planes. More formally: For positive integers $t \leq k
\leq v$ and $\lambda$, we define a \mbox{\emph{$t$-$(v,k,\lambda)$
design}} to be a finite incidence structure
\mbox{$\mathcal{D}=(X,\mathcal{B},I)$}, where $X$ denotes a set of
\emph{points}, $\left| X \right| =v$, and $\mathcal{B}$ a set of
\emph{blocks}, $\left| \mathcal{B} \right| =b$, with the following
regularity properties: each block $B \in \mathcal{B}$ is incident
with $k$ points, and each \mbox{$t$-subset} of $X$ is incident with
$\lambda$ blocks. A \emph{flag} of $\mathcal{D}$ is an incident
point-block pair $(x,B) \in I$ with $x \in X$ and $B \in
\mathcal{B}$.

For historical reasons, a \mbox{$t$-$(v,k,\lambda)$ design} with
$\lambda =1$ is called a \emph{Steiner \mbox{$t$-design}} (sometimes
also a \emph{Steiner system}). We note that in this case each block
is determined by the set of points which are incident with it, and
thus can be identified with a \mbox{$k$-subset} of $X$ in a unique
way. If $t<k<v$, then we speak of a \emph{non-trivial} Steiner
\mbox{$t$-design}. As a simple example, the vector space
$\mathbb{Z}_2^n$ ($n \geq 3$) with block set $\mathcal{B}$ taken to
be the set of all subsets of four distinct elements of
$\mathbb{Z}_2^n$ whose vector sum is zero is a (boolean) Steiner
\mbox{$3$-$(2^n,4,1)$ design}. There are many infinite classes of
Steiner \mbox{$t$-designs} for $t=2$ and $3$, however for $t=4$ and
$5$ only a finite number are known. For a detailed treatment of
combinatorial designs, we refer
to~\cite{bjl99,crc06,hall86,hupi85,stin04}. In
particular,~\cite{bjl99,crc06} provide encyclopedic accounts of key
results and contain existence tables with known parameter sets.

In what follows, we are interested in \mbox{$t$-designs} which admit
groups of automorphisms with sufficiently strong symmetry properties
such as transitivity on the blocks or on the flags. We consider
automorphisms of a \mbox{$t$-design} $\mathcal{D}$ as pairs of
permutations on $X$ and $\mathcal{B}$ which preserve incidence, and
call a group \mbox{$G \leq \mbox{Aut} (\mathcal{D})$} of
automorphisms of $\mathcal{D}$ \emph{block-transitive} (respectively
\emph{flag-transitive}, \emph{point \mbox{$t$-transitive}},
\emph{point \mbox{$t$-homogeneous}}) if $G$ acts transitively on the
blocks (respectively transitively on the flags,
\mbox{$t$-transitively} on the points, \mbox{$t$-homogeneously} on
the points) of $\mathcal{D}$. For short, $\mathcal{D}$ is said to
be, e.g., block-transitive if $\mathcal{D}$ admits a
block-transitive group of automorphisms.


\section{Basic Properties and Existence Results}\label{basics}

We give some basic properties and known results concerning the
existence of \mbox{$t$-designs} which are important for the
remainder of the paper.

\smallskip

If $\mathcal{D}=(X,\mathcal{B},I)$ is a \mbox{$t$-$(v,k,\lambda)$}
design with $t \geq 2$, and $x \in X$ arbitrary, then the
\emph{derived} design with respect to $x$ is
\mbox{$\mathcal{D}_x=(X_x,\mathcal{B}_x, I_x)$}, where $X_x = X
\backslash \{x\}$, \mbox{$\mathcal{B}_x=\{B \in \mathcal{B}:
(x,B)\in I\}$} and $I_x= I \! \mid _{X_x \times \, \mathcal{B}_x}$.
In this case, $\mathcal{D}$ is also called an \emph{extension} of
$\mathcal{D}_x$. Obviously, $\mathcal{D}_x$ is a
\mbox{$(t-1)$-$(v-1,k-1,\lambda)$} design.

For \mbox{$\mathcal{D}=(X,\mathcal{B},I)$} a Steiner
\mbox{$t$-design} with \mbox{$G \leq \mbox{Aut} (\mathcal{D})$}, let
$G_x$ denote the stabilizer of a point $x \in X$, and $G_B$ the
setwise stabilizer of a block $B \in \mathcal{B}$. For $x, y \in X$
and $B \in \mathcal{B}$, we define $G_{xy}= G_x \cap G_y$ and
$G_{xB}= G_x \cap G_B$.

For any $x \in \mathbb{R}$, let $\lfloor x \rfloor$ denote the
greatest positive integer which is at most $x$.

All other notations remain as defined in Section~\ref{designs}.

\smallskip

Basic necessary conditions for the existence of \mbox{$t$-designs}
can be obtained via elementary counting arguments (see, for
instance,~\cite{bjl99}):

\begin{prop}\label{s-design}
Let $\mathcal{D}=(X,\mathcal{B},I)$ be a \mbox{$t$-$(v,k,\lambda)$}
design, and for a positive integer $s \leq t$, let $S \subseteq X$
with $\left|S\right|=s$. Then the total number of blocks incident
with each element of $S$ is given by
\[\lambda_s = \lambda \frac{{v-s \choose t-s}}{{k-s \choose t-s}}.\]
In particular, for $t\geq 2$, a \mbox{$t$-$(v,k,\lambda)$} design is
also an \mbox{$s$-$(v,k,\lambda_s)$} design.
\end{prop}

\noindent It is customary to set $r:= \lambda_1$ denoting the total
number of blocks incident with a given point (referring to the
`replication number' from statistical design of experiments, one of
the origins of design theory).

\begin{Cor} \label{Comb_t=5}
Let $\mathcal{D}=(X,\mathcal{B},I)$ be a \mbox{$t$-$(v,k,\lambda)$}
design. Then the following holds:
\begin{enumerate}

\item[{(a)}] $bk = vr.$

\smallskip

\item[{(b)}] $\displaystyle{{v \choose t} \lambda = b {k \choose t}.}$

\smallskip

\item[{(c)}] $r(k-1)=\lambda_2(v-1)$ for $t \geq 2$.

\end{enumerate}
\end{Cor}

\begin{Cor}\label{divCond}
Let $\mathcal{D}=(X,\mathcal{B},I)$ be a \mbox{$t$-$(v,k,\lambda)$}
design. Then
\[\lambda {v-s \choose t-s} \equiv \, 0\; \emph{(mod}\;\, {k-s \choose t-s})\]
for each positive integer $s \leq t$.
\end{Cor}

\smallskip

For non-trivial Steiner \mbox{$t$-designs} lower bounds for $v$ in
terms of $k$ and $t$ can be given (see
P.~Cameron~\cite[Thm.\,3A.4]{Cam1976}, and
J.~Tits~\cite[Prop.\,2.2]{Tits1964}):

\begin{thm}
\label{Cam} If $\mathcal{D}=(X,\mathcal{B},I)$ is a non-trivial
Steiner \mbox{$t$-design}, then the following holds:
\begin{enumerate}

\item[{(a)}] \emph{(Tits~1964):} \hspace{0.2cm} $v\geq (t+1)(k-t+1).$

\smallskip

\item[{(b)}] \emph{(Cameron~1976):} \hspace{0.2cm} \mbox{$v-t+1 \geq (k-t+2)(k-t+1)$} for $t>2$. If equality
holds, then
\smallskip
$(t,k,v)=(3,4,8),(3,6,22),(3,12,112),(4,7,23)$, or $(5,8,24)$.
\end{enumerate}
\end{thm}

\noindent We note that (a) is stronger for $k<2(t-1)$, while (b) is
stronger for \linebreak \mbox{$k>2(t-1)$}. For $k=2(t-1)$ both
assert that $v \geq t^2-1$.

\smallskip

The following result by R.~Fisher~\cite{fish40} is classical,
generally known as ``Fisher's Inequality''\index{Fisher's
Inequality}:

\begin{thm}{\em (Fisher~1940).}\label{FisherIn}
If $\mathcal{D}=(X,\mathcal{B},I)$ is a non-trivial
\mbox{$2$-$(v,k,\lambda)$} design, then $b \geq v$, that is, there
are at least as many blocks as points in $\mathcal{D}$.
\end{thm}

An important generalization to arbitrary \mbox{$t$-designs} is due
to D.~Ray-Chaudhuri and R.~Wilson~\cite[Thm.\,1]{Ray-ChWil1975}:

\begin{thm}{\em (Ray-Chaudhuri \& Wilson~1975).}\label{RayCh}
Let $\mathcal{D}=(X,\mathcal{B},I)$ be a \mbox{$t$-$(v,k,\lambda)$}
design. If $t$ is even, say $t=2s$, and $v \geq k+s$, then $b \geq
{v\choose s}$. If $t$ is odd, say $t = 2s+1$, and $v-1 \geq k+s$,
then $b \geq 2{v-1\choose s}$.
\end{thm}

Exploration of the construction of \mbox{$t$-designs} for large
values of $t$ lead to L.~Teirlinck's celebrated
theorem~\cite{Teir1987}, one of the major results in design theory:

\begin{thm}{\em (Teirlinck~1987).}\label{Teirl}
For every positive integer value of $t$, there exists a non-trivial
\mbox{$t$-design}.
\end{thm}

However, although Teirlinck's recursive methods are constructive,
they only produce examples with tremendously large values of
$\lambda$. Until now no non-trivial Steiner \mbox{$t$-design} with
$t>5$ has been constructed.

\begin{problem*}
Does there exist any non-trivial Steiner \mbox{$6$-design}?
\end{problem*}


\section{Approach via Symmetry}

Besides recursive and set-theoretical approaches, many existence
results for \mbox{$t$-designs} with large $t$ have been obtained in
recent years by the method of orbiting under a group (see,
e.g.,~\cite{betten2004},~\cite[§II.4]{crc06}). Specifically, the
consideration of \mbox{$t$-designs} which admit groups of
automorphisms with sufficiently strong symmetry properties seems to
be of great importance in our context - quoting arguably Gian-Carlo
Rota~\cite{Rota1964}:

\begin{quote}
``A combinatorial object without symmetries doesn't exist
 - by definition.''
\end{quote}

\smallskip

We first state (cf.~\cite{Tuan2002}):
\begin{prop}
Let $t$ be a positive integer, and $G$ a finite (abstract) group.
Then there is a \mbox{$t$-design} such that the full group
$\emph{Aut}(\mathcal{D})$ of automorphisms has a subgroup isomorphic
to $G$.
\end{prop}

One of the early important results regarding highly symmetric
designs is due to R.~Block~\cite[Thm.\,2]{Block1965}:

\begin{prop}{\em (Block~1965).}\label{BlocksLemma2}
Let $\mathcal{D}=(X,\mathcal{B},I)$ be a non-trivial
\mbox{$t$-$(v,k,\lambda)$} design with $t \geq 2$. If $G \leq
\emph{Aut}(\mathcal{D})$ acts block-transitively on $\mathcal{D}$,
then $G$ acts point-transitively on $\mathcal{D}$.
\end{prop}

For a \mbox{$2$-$(v,k,1)$} design $\mathcal{D}$, it is elementary
that the point \mbox{$2$-transitivity} of \mbox{$G \leq
\mbox{Aut}(\mathcal{D})$} implies its flag-transitivity. For
\mbox{$2$-$(v,k,\lambda)$} designs, this implication remains true if
$r$ and $\lambda$ are relatively prime (see,
e.g.,~\cite[Chap.\,2.3,\,Lemma\,8]{demb68}). However, for
\mbox{$t$-$(v,k,\lambda)$} designs with $t \geq 3$, it can be
deduced from Proposition~\ref{BlocksLemma2} that always the converse
holds (see~\cite{Buek1968} or~\cite[Lemma\,2]{Hu2001}):

\begin{prop}\label{flag2trs}
Let $\mathcal{D}=(X,\mathcal{B},I)$ be a non-trivial
\mbox{$t$-$(v,k,\lambda)$} design with $t \geq 3$. If \mbox{$G \leq
\emph{Aut}(\mathcal{D})$} acts flag-transitively on $\mathcal{D}$,
then $G$ acts point \linebreak \mbox{$2$-transitively} on
$\mathcal{D}$.
\end{prop}

Investigating highly symmetric \mbox{$t$-designs} for large values
of $t$, P.~Cameron and C.~Praeger~\cite[Thm.\,2.1]{CamPrae1993}
deduced from Theorem~\ref{RayCh} and Proposition~\ref{BlocksLemma2}
the following assertion:

\begin{prop}{\em (Cameron \& Praeger~1993).}\label{flag3hom}
Let $\mathcal{D}=(X,\mathcal{B},I)$ be a \mbox{$t$-$(v,k,\lambda)$}
design with $t\geq 2$. Then, the following holds:

\begin{enumerate}

\item[{(a)}] If \mbox{$G \leq \emph{Aut}(\mathcal{D})$} acts block-transitively on $\mathcal{D}$,
then $G$ also acts point \linebreak \mbox{$\lfloor t/2
\rfloor$-homogeneously} on $\mathcal{D}$.

\smallskip

\item[{(b)}] If \mbox{$G \leq \emph{Aut}(\mathcal{D})$} acts flag-transitively on $\mathcal{D}$,
then $G$ also acts point \linebreak \mbox{$\lfloor (t+1)/2
\rfloor$-homogeneously} on $\mathcal{D}$.

\end{enumerate}
\end{prop}

As for $t \geq 7$ the flag-transitivity, respectively for $t \geq 8$
the block-transitivity of \mbox{$G \leq \mbox{Aut} (\mathcal{D})$}
implies at least its point \mbox{$4$-homogeneity}, they obtained the
following restrictions as a fairly direct consequence of the finite
simple group classification (cf.~\cite[Thm.\,1.1]{CamPrae1993}):

\begin{thm}{\em (Cameron \& Praeger~1993).}
Let $\mathcal{D}=(X,\mathcal{B},I)$ be a \mbox{$t$-$(v,k,\lambda)$}
design. If \mbox{$G \leq \emph{Aut}(\mathcal{D})$} acts
block-transitively on $\mathcal{D}$ then $t \leq 7$, while if
\mbox{$G \leq \emph{Aut}(\mathcal{D})$} acts flag-transitively on
$\mathcal{D}$ then $t \leq 6$.
\end{thm}

Moreover, they formulated the following far-reaching conjecture
(cf.~\cite[Conj.\,1.2]{CamPrae1993}):

\begin{conj}{\em (Cameron \& Praeger~1993).}
There are no non-trivial block-transitive \mbox{$6$-designs}.
\end{conj}

It has been shown recently by the author~\cite{hu07} that no
non-trivial flag-transitive Steiner \mbox{$6$-design} can exist.
Moreover, he classified all flag-transitive Steiner
\mbox{$t$-designs} with $t>2$
(see~\cite{Hu2001,Hu2005,Hu_Habil2005,hu07,hu07a} and~\cite{Hu2008}
for a monograph). These results make use of the classification of
all finite \mbox{$2$-transitive} permutation groups, which itself
relies on the finite simple group classification. The
characterizations answer a series of 40-year-old problems and
generalize theorems of J.~Tits~\cite{Tits1964} and
H.~Lüneburg~\cite{Luene1965}. Earlier, F.~Buekenhout,
A.~Delandtsheer, J.~Doyen, P.~Kleidman, M.~Liebeck, and
J.~Saxl~\cite{Buek1990,Del2001,Kleid1990,Lieb1998,Saxl2002} had
essentially characterized all finite flag-transitive linear spaces,
that is flag-transitive Steiner \mbox{$2$-designs}. Their result,
which also relies on the finite simple group classification, starts
with a classical result of Higman and McLaughlin~\cite{HigMcL1961}
and uses the O'Nan-Scott Theorem for finite primitive permutation
groups. For the incomplete case with a \mbox{$1$-dimensional} affine
group of automorphisms, we refer to~\cite[Sect.\,4]{Buek1990}
and~\cite[Sect.\,3]{Kant1993}.


\section{Non-Existence of Block-transitive Steiner \mbox{$6$-Designs}}\label{block6designs}

We assert the following main result:

\begin{mthm*}\label{mainthm}
Let $\mathcal{D}=(X,\mathcal{B},I)$ be a non-trivial Steiner
\mbox{$6$-design}. Then \mbox{$G \leq \emph{Aut}(\mathcal{D})$}
cannot act block-transitively on $\mathcal{D}$, except possibly when
$G= P \mathit{\Gamma} L(2,p^e)$ with $p=2$ or $3$ and $e$ is an odd
prime power.
\end{mthm*}

We will briefly outline the main ingredients of the proof. The long
and technical details will appear elsewhere:

\smallskip

$\bullet \;$ In order to investigate block-transitive Steiner
\mbox{$6$-designs}, we can in view of Proposition~\ref{flag3hom}~(a)
make use of the classification of all finite \mbox{$3$-homogeneous}
permutation groups, which itself relies on the finite simple group
classification
(cf.~\cite{Cam1981,Gor1982,Kant1972,Lieb1987,LivWag1965}). The list
of groups which have to be examined is as follows:

Let $G$ be a finite \mbox{$3$-homogeneous} permutation group on a
set $X$ with $\left|X\right| \geq 4$. Then $G$ is either of

{\bf (A) Affine Type:} $G$ contains a regular normal subgroup $T$
which is elementary Abelian of order $v=2^d$. If we identify $G$
with a group of affine transformations
\[x \mapsto x^g+u\]
of $V=V(d,2)$, where $g \in G_0$ and $u \in V$, then particularly
one of the following occurs:

\begin{enumerate}

\smallskip

\item[(1)] $G \cong AGL(1,8)$, $A \mathit{\Gamma} L(1,8)$, or $A \mathit{\Gamma} L(1,32)$

\smallskip

\item[(2)] $G_0 \cong SL(d,2)$, $d \geq 2$

\smallskip

\item[(3)] $G_0 \cong A_7$, $v=2^4$

\end{enumerate}

\smallskip

or

\smallskip

{\bf (B) Almost Simple Type:} $G$ contains a simple normal subgroup
$N$, and \mbox{$N \leq G \leq \mbox{Aut}(N)$}. In particular, one of
the following holds, where $N$ and $v=|X|$ are given as follows:
\begin{enumerate}

\smallskip

\item[(1)] $A_v$, $v \geq 5$

\smallskip

\item[(2)] $PSL(2,q)$, $q>3$, $v=q+1$

\smallskip

\item[(3)] $M_v$, $v=11,12,22,23,24$ \hfill (Mathieu groups)

\smallskip

\item[(4)] $M_{11}$, $v=12$

\end{enumerate}

\smallskip

We note that if $q$ is odd, then $PSL(2,q)$ is
\mbox{$3$-homogeneous} for \mbox{$q \equiv 3$ (mod $4$)}, but not
for \mbox{$q \equiv 1$ (mod $4$)}, and hence not every group $G$ of
almost simple type satisfying (2) is \mbox{$3$-homogeneous} on $X$.

\smallskip

$\bullet \;$ If \mbox{$G \leq \mbox{Aut}(\mathcal{D})$} acts
block-transitively on any Steiner \mbox{$t$-design} $\mathcal{D}$
with $t\geq 6$, then in particular $G$ acts point
\mbox{$2$-transitively} on $\mathcal{D}$ by
Proposition~\ref{flag3hom}~(a). Applying Lemma~\ref{Comb_t=5}~(b)
yields the equation
\[b=\frac{{v \choose t}}{{k \choose
t}}=\frac{v(v-1) \left|G_{xy}\right|}{\left| G_B \right|},\] where
$x$ and $y$ are two distinct points in $X$ and $B$ is a block in
$\mathcal{B}$. Combined with the combinatorial properties in
Section~\ref{basics}, this arithmetical condition yields in some of
the cases under consideration immediately strong results. In other
cases, particular Diophantine equations arise which have to be
examined in more detail.

\smallskip

$\bullet \;$ As for the flag-transitive treatment
(cf.~\cite{hu07,hu07a}), the projective group containing $PSL(2,q)$
-- although group-theoretically well understood -- requires some
complicated analysis in this context. This includes a detailed
consideration of the orbit-lengths from the action of subgroups of
$PSL(2,q)$ on the points of the projective line. The cases excluded
from the theorem remain elusive. However, it seems to be very
unlikely that admissible parameter sets of Steiner
\mbox{$6$-designs} can be found in view of the arithmetical
conditions that are imposed in these cases.

\bibliographystyle{amsplain}
\bibliography{XbibCamPraeg}
\end{document}